\documentclass{ACmod}

\DeclareFontFamily{U}{rsf}{}
\DeclareFontShape{U}{rsf}{m}{n}{
  <5> <6> rsfs5 <7> <8> <9> rsfs7 <10-> rsfs10}{}
\DeclareMathAlphabet{\mathscr}{U}{rsf}{m}{n}

\makeatletter
\def\operator@font{\sf}
\makeatother

\newcommand{\g}{\mathfrak{g}}

\newcommand{\sH}{{\mathcal H}}

\newcommand{\sX}{{\mathfrak X}}

\newcommand{\cE}{{\mathscr E}}
\newcommand{\cF}{{\mathscr F}}

\newcommand{\cO}{{\mathscr O}}

\newcommand{\bE}{{\mathbb{E}}}

\newcommand{\tW}{{\widetilde{W}}}
\newcommand{\tXg}{{\widetilde{X^g}}}
\newcommand{\tXh}{{\widetilde{X^h}}}
\newcommand{\tXk}{{\widetilde{X^k}}}
\newcommand{\tXgh}{{\widetilde{X^{gh}}}}
\newcommand{\tXgch}{{\widetilde{X^{g,h}}}}
\newcommand{\tX}{\widetilde{X}}

\newcommand{\orb}{{\mathsf{orb}}}

\newcommand{\HKR}{{\mathsf{HKR}}}

\newcommand{\chk}{{\scriptscriptstyle\vee}}

\newcommand{\bT}{\mathbb{T}}
\newcommand{\bD}{\mathbb{D}}
\newcommand{\bN}{\mathbb{N}}

\newcommand{\bL}{{\mathbb{L}}}

\DeclareMathOperator{\Sym}{Sym}

\DeclareMathOperator{\Cone}{Cone}

\DeclareMathOperator{\HH}{HH}
\DeclareMathOperator{\HT}{HT}

\DeclareMathOperator{\Td}{td}

\DeclareMathOperator{\ord}{ord}

\DeclareMathOperator{\id}{id}

\DeclareMathOperator{\rank}{rk}
\DeclareMathOperator{\Ext}{Ext}

\DeclareMathOperator{\Hom}{Hom}

\DeclareMathOperator{\codim}{codim}

\newcommand{\ra}{\rightarrow}

\newcommand{\lra}{\longrightarrow}

\newcommand{\C}{\mathbb{C}}

\newcommand{\Z}{\mathbb{Z}}

\newcommand{\iso}{\cong}

\renewcommand{\phi}{\varphi}

\usepackage{tikz-cd}

\usepackage{xypic}
\usepackage{tikz}
\usepackage[tocflat]{tocstyle}
\usetikzlibrary{positioning}
\usepackage{amscd}
\usetikzlibrary{arrows,shapes,positioning}
\usetikzlibrary{decorations.markings}
\tikzstyle arrowstyle=[scale=1]
\tikzstyle directed=[postaction={decorate,decoration={markings,
    mark=at position .65 with {\arrow[arrowstyle]{stealth}}}}]
\tikzstyle reverse directed=[postaction={decorate,decoration={markings,
    mark=at position .65 with {\arrowreversed[arrowstyle]{stealth};}}}]
\usepackage[tocflat]{tocstyle}
\tikzset{cross/.style={cross out, draw=black, minimum size=2*(#1-\pgflinewidth), inner sep=0pt, outer sep=0pt},
cross/.default={1pt}}

\begin{document}

\title{The cup product in orbifold Hochschild cohomology}

\author[C\u ald\u araru and Huang]{%
Andrei C\u ald\u araru and Shengyuan Huang}

\address{
\textsc{Andrei C\u ald\u araru: Mathematics Department,
University of Wisconsin--Madison
\newline
480 Lincoln Drive, Madison, WI 53706--1388, USA}
\newline
Email address: \texttt{andreic@math.wisc.edu}
\newline
\newline
\textsc{Shengyuan Huang: Mathematics Department,
University of Wisconsin--Madison
\newline
480 Lincoln Drive, Madison, WI 53706--1388, USA}
\newline
Email address: \texttt{shuang279@wisc.edu}}

\begin{abstract}
{\sc Abstract:}
We study the multiplicative structure of orbifold Hochschild
cohomology in an attempt to generalize the results of
Kontsevich and Calaque-Van den Bergh relating the Hochschild and
polyvector field cohomology rings of a smooth variety.

We introduce the concept of linearized derived scheme, and we
argue that when $X$ is a smooth algebraic variety and $G$ is a finite
abelian group acting on $X$, the derived fixed locus $\widetilde{X}^G$
admits an HKR linearization.  This allows us to define a product on the
cohomology of polyvector fields of the orbifold $[X/G]$.  We analyze
the obstructions to associativity of this product and show that they
vanish in certain special cases. We conjecture that in these cases the
resulting polyvector field cohomology ring is isomorphic to the
Hochschild cohomology of $[X/G]$.

Inspired by mirror symmetry we introduce a bigrading on the Hochschild
homology of Calabi-Yau orbifolds.  We propose a conjectural product
which respects this bigrading and simplifies the previously introduced
product.
\end{abstract}

\maketitle
\setcounter{tocdepth}{1}
\tableofcontents
\vfill

\section{Introduction}
\paragraph
\label{1.1}
Let $X$ be a smooth algebraic variety over a field of characteristic
zero.  The HKR map is a graded vector space isomorphism
\[ \HT^*(X) \stackrel{\sim}{\lra} \HH^*(X) \]
between the  {\em polyvector field cohomology} $\HT^*(X)$ of $X$,
\[ \HT^*(X)=\bigoplus_{p+q=*} H^p(X,\wedge^qT_X), \]
and the Hochschild cohomology $\HH^*(X)$ of $X$.

Both spaces above are graded commutative rings: polyvector field
cohomology classes can be multiplied using the wedge product on
$\wedge^* TX$ and cup product on cohomology, while Hochschild
cohomology classes can be composed using the Yoneda product. However,
the HKR isomorphism is not a ring map in general.

Kontsevich~\cite{K} claimed that the rings $\HT^*(X)$ and $\HH^*(X)$
{\em are} in fact isomorphic, via a modification of the HKR
isomorphism. This result was later proved by Calaque and Van den
Bergh~\cite{CV}.

\paragraph
The ring $\HT^*(X)$ is bigraded, and the product respects this
bigrading. Moreover, the Hochschild cochain complex carries a
filtration given by order of polydifferential operators, which in turn
induces a filtration on Hochschild cohomology. Kontsevich's claim
(along with the explicit formula for the corrected ring isomorphism)
can be interpreted as saying that this filtration admits a
multiplicative splitting, yielding a bigrading on $\HH^*(X)$ which
refines the usual grading.

\paragraph
The problem of understanding an analogue of Kontsevich's claim for
orbifolds has been open for at least 20 years.  The most recent
(negative) progress is due to
Negron-Schedler~\cite{NS} who argue that the
Hochschild cochain complex of an orbifold does not satisfy a formality
result similar to the one Kontsevich used.  However, this does not
rule out the possibility that an analogue of the Kontsevich claim
holds for a {\em corrected} filtration from the one they study.

This problem is particularly interesting in view of its connections to
Ruan's crepant resolution conjecture. For example, getting a good
understanding of the orbifold Hochschild cohomology product would
explain the matching between the cohomology ring of the Hilbert scheme
of $n$-points on a K3 surface $S$ and the Chen-Ruan orbifold
cohomology ring of $[S^n/\Sigma_n]$, as observed by
Fantechi-G\"ottsche~\cite{FG}.

\paragraph
Consider a global quotient orbifold $[X/G]$, where $X$ is a smooth
algebraic variety and $G$ is a finite group acting on $X$.  Arinkin,
C\u{a}ld\u{a}raru, and Hablicsek~\cite{ACH} gave an explicit
decomposition of the Hochschild cohomology of $[X/G]$ as a graded
vector space
\[
\HH^*([X/G]) \stackrel{\sim}{\longleftarrow} \left (\bigoplus_{g\in
    G}\bigoplus_{p+q=*} H^{p-c_g}(X^g,\wedge^q
  T_{X^g}\otimes\omega_g)\right )^G\!\!\!\!\!,
\]
where $X^g$ is the fixed locus of $g\in G$, $c_g$ is the codimension
of $X^g$ in $X$, and $\omega_g$ is the dualizing sheaf of the
inclusion $X^g\hookrightarrow X$. The right hand side is the natural
analogue of polyvector field cohomology for orbifolds, and the above
map can be regarded as the natural generalization of the HKR
isomorphisms for global quotient orbifolds.

We define
\[ \HT^*(X;G)=\left (\bigoplus_{g\in G}\bigoplus_{p+q=*}
    H^{p-c_g}(X^g,\wedge^q T_{X^g}\otimes\omega_g)\right ). \]
Note that $\HT^*(X;G)$ carries a natural $G$ action, and we set
\[  \HT^*([X/G])=\HT^*(X;G)^G.\]

\paragraph
As stated above, for a smooth variety $X$ there is an obvious
associative product on $\HT^*(X)$. However, when $G$ is non-trivial,
it is not at all obvious what the analogous product structure should
be on $\HT^*([X/G])$. Understanding candidates for such a product is
the goal of this paper.

We will define two operations on $\HT^*(X;G)$. The first one is
inspired closely by the construction of the product of distributions
on a Lie algebra. We can prove the associativity of this operation
when certain cohomology classes vanish. In particular, this product
will be associative when $X$ is affine. However, in general we can
only define this operation for abelian $G$, and we are not able to
construct a multiplicative bigrading. A special case of this
construction, which works with a non-abelian group $G$, will be
discussed in Section 7.

The second definition mimics a construction of Fantechi-G\"{o}ttsche,
and the resulting formulas conjecturally simplify the first
construction above. The operation we construct respects a natural
bigrading, and it behaves well in all the examples we can
compute. However, we cannot prove a general criterion for
associativity of this operation.

\paragraph{\bf Theorem A.} {\em
Suppose $[X/G]$ is a global quotient orbifold, where $X$ is a smooth
algebraic variety and $G$ is a finite {\em abelian} group acting on
$X$. Then the construction in Section 3 defines an operation on
$\HT^*(X;G)$ which recovers the wedge product on $\HT^*(X)$ when $G$
is trivial.

This operation is associative if the Bass-Quillen class~(\ref{def:bq})
associated to the sequence of closed embeddings
$X^{g,h}\hookrightarrow X^{g}\hookrightarrow X$ vanishes for all
$g,h\in G$. (The abstract Bass-Quillen class was introduced by the
second author in~\cite{H}.)}

In particular, the product is associative when $X$ is affine, or when
$X$ is an abelian variety and $G = \mathbb{Z}/2\mathbb{Z}$ acting by
negation.

\paragraph{\bf Conjecture A.}
\label{conj}
{\em
Let $[X/G]$ be a global quotient orbifold with $G$ abelian.
Then the HKR isomorphism can be corrected to an isomorphism of rings
\[ \HH^*([X/G])\cong \HT^*([X/G]).\]
We will give in Section 7 an explicit conjectural formula for the corrected
isomorphism.}

\paragraph
The relationship to bigradings can be understood from the point of
view of mirror symmetry. Recall that under mirror symmetry the
Hochschild cohomology of a smooth Calabi-Yau variety $X$ corresponds
to the singular cohomology of its mirror $\check{X}$.  Since the
latter admits a multiplicative bigrading (given by the Hodge
decomposition) it is natural to expect that $\HH^*(X)$ should also
carry a multiplicative bigrading.

For orbifolds, Chen and Ruan~\cite{CR} used ideas from quantum
cohomology to define an analogue of singular cohomology for
orbifolds.  The Chen-Ruan orbifold cohomology space is defined to be
\[ H^*([X/G],\mathbb{C})=\left(\bigoplus_{g\in G}
  H^{*-2\iota(g)}(X^g,\mathbb{C})\right)^G, \]
where $\iota(g)\in\mathbb{Q}$ is the age of an element $g\in{G}$. This
space is endowed with a natural, graded commutative product. Fantechi
and G\"{o}ttsche~\cite{FG} later gave an explicit description of this product
and noted that it has a natural bigrading induced by the Hodge
decompositions on the fixed loci $X^g$:
\[ \left(\bigoplus_{g\in G}H^{*-2\iota(g)}(X^g,\mathbb{C})\right)^G=
  \left(\bigoplus_{g\in G}\bigoplus_{p,q}H^{p-\iota(g)}(X^g,
  \wedge^{q-\iota(g)}\Omega_{X^g})\right)^G.\]
The mirror symmetry intuition above suggests that the Hochschild
cohomology of Calabi-Yau global quotient orbifolds should also carry a
multiplicative bigrading.

\paragraph
We define a new bigrading on the cohomology of orbifold polyvector
fields by setting
\[ \HT^{p,q}(X;G)=\bigoplus_{g\in G}
  H^{p-\iota(g)}(X^g,\wedge^{q+\iota(g)-c_g} T_{X^g}\otimes\omega_g).
\]
Based on examples we compute in this paper
we propose a conjectural way to simplify the product in Theorem A for
Calabi-Yau global quotient orbifolds.  This simplified product will
preserve the above bigrading.

\paragraph {\bf Plan of the paper.}
In Section 2 we briefly summarize Kontsevich's Theorem for smooth
algebraic varieties. Then we turn to the global quotient orbifold
case, and recall the identification of orbifold polyvector fields with
orbifold Hochschild cohomology via the HKR isomorphism.

Section 3 is devoted to the definition of our product structure on
orbifold polyvector fields.

Section 4 contains the technical details based on derived intersection
and formality of derived schemes. We need to use them to make the
definition work. The definition of our product works at least when $G$
is abelian, and also works for examples of Hilbert schemes of $n$
points as will be described in Section 7.

Section 5 presents the strategy for proving the associativity of our
product when the Bass-Quillen class vanishes. We explain how to reduce
the question of associativity to two statements, Propositions
\ref{prop assoc 1} and \ref{prop assoc 2}, which are proved in Section 6.

In Section 7 we propose an ansatz which allows us to simplify the
formulas for our product. We compute several examples of the
simplified product, leading us to a definition of the bigrading
described above on the cohomology of orbifold polyvector fields. The
simplified product preserves this bigrading. At the end of this
paper, we list several remaining problems directly related to
Conjecture A.

\paragraph {\bf Conventions.} All the algebraic varieties in this
paper are smooth over a field of characteristic zero.

\paragraph {\bf Acknowledgments.} We would like to thank Dima Arinkin
for patiently listening to the various problems we ran into at
different stages of the project, and for providing insight.

The authors were partially supported by the National Science
Foundation through grant number DMS-1811925.

\section{Background}
In this section we review the theorem of Kontsevich and Calaque-Van
den Bergh for smooth algebraic varieties, and the construction of the
orbifold HKR isomorphism. A few new definitions and interpretations in
this section will be important throughout the paper.

\paragraph{\bf The HKR isomorphism.}
Let $X$ be a smooth algebraic variety. There is an HKR isomorphism
\[ \HH^*(X)\cong\bigoplus_{p+q=*}H^{p}(X,\wedge^q T_X) \]
that identifies the Hochschild cohomology of $X$ with the cohomology
of polyvector fields as graded vector spaces. More precisely, we have
a sheaf-level HKR isomorphism in the derived category of $X$
\[ \Delta^*\Delta_*\cO_X\cong\Sym_{\cO_X}\Omega_X[1], \]
where $\Delta: X\hookrightarrow X\times X$ is the diagonal embedding.
We get the desired isomorphism on cohomology by applying
$\Hom(\,-\,,\cO_X)$ to this isomorphism of sheaves.

\paragraph
The HKR isomorphism has an interpretation in the language of derived
schemes. The complex $\Delta^*\Delta_*\cO_X$ admits a graded
commutative product, which makes it into the structure complex of a
derived scheme, the derived self-intersection $X\times^R_{X\times X}
X$.  This is also known as $LX$, the derived loop space of $X$.

On the other hand $\Sym \Omega_X[1]$ is the structure complex
of the shifted tangent bundle $\bT_{X}[-1]$.  We can then restate the
HKR isomorphism as an isomorphism of derived schemes over $X$,
\[ \xymatrix{ \bT_X[-1]\ar[r]^{\cong~~~~~~~~} &
    LX=X\times^R_{X\times X} X. } \]

\paragraph
It was observed by Kapranov and Kontsevich that there is a Lie
theoretic interpretation of the HKR isomorphism. The derived loop
space $LX$ has the structure of a derived group
scheme over $X$, and the relative normal bundle $N_{X/LX}=\bT_X[-1]$ is
its Lie algebra~\cite{Ka}.  The HKR isomorphism can be thought of as a
version of the exponential map $\bT_X[-1]\rightarrow LX$
\cite{CalRo}, where $\bT_X[-1]$ is the total space of the
shifted tangent bundle.

The exponential map relates the derived Lie group $LX$ with its Lie
algebra $N_{X/LX}=\bT_X[-1]$. In general, suppose $\tX$ is an
arbitrary derived scheme which is not necessary a derived group. We
can still consider the total space $\bN_{X/\tX}$ of
the normal bundle of $X$ in $\tX$, where $X$ is the
underlying classical scheme $X\hookrightarrow\tX$. This
leads to the following definition.

\begin{Definition}
For a derived scheme $\tX$, the {\em linearization} $\bL_{\tX}$ of $\tX$
is defined to be the total space of the normal bundle $N_{X/\tX}$,
where $X$ is the underlying classical scheme $X\hookrightarrow\tX$.  A
choice of isomorphism $\bL_{\tX}\cong\tX$ (if one exists) will be
called a {\em linearization} of $\tX$. 
\end{Definition}
\bigskip

For example, consider $X\hookrightarrow LX$. Then
\[ \bL_{LX}=\bN_{X/LX}=\bT_X[-1], \]
and the HKR isomorphism provides a linearization of $LX$.

\paragraph
We need to address a technical detail about the above isomorphism.
The linearization $\bL_{\tX}$ of a derived scheme $\tX$ is by
definition the total space of the normal bundle $N_{X/\tX}$, hence it
comes with a natural projection which makes it a scheme over $X$.
Moreover, this projection splits the inclusion $X\hookrightarrow \tX$.
However, in general $\tX$ may not admit such a projection.  This
explains why in general here is no way to define an isomorphism
$\bL_{\tX}\cong \tX$ over $X$.  If we hope to define an
isomorphism $\bL_{\tX}\cong \tX$, we usually need to find a natural
base scheme $Y$ such that both $\tX$ and $\bL_{\tX}$ are affine over
$Y$. Then we can consider the structure complex of $\bL_{\tX}$ and
$\tX$ as $\cO_Y$-algebras. There is a bijection between the set of
isomorphisms $\cO_{\bL_{\tX}}\cong\cO_{\tX}$ in the derived category
of $Y$ and the set of isomorphisms $\bL_{\tX}\cong\tX$ over $Y$. A
choice of such an isomorphism will be called a linearization of $\tX$
over $Y$.

The HKR isomorphism
\[ \xymatrix{
  \exp:\bT_{X}[-1]=\bL_{LX}\ar[r]^{~~~~~~~~~~~~\cong} & LX } \]
linearizes the derived loop space $LX$ over $X$. Here
$LX=X\times^R_{X\times X}X$ is to be viewed as a scheme over $X$ via
one of the two projection maps onto the left or right factors
$X\times^R_{X\times X}X\rightarrow X$.

\paragraph
\label{formal}
For most derived schemes $\tX$, even if they are affine over a scheme
$Y$, it is not true that $\tX$ is isomorphic to $\bL_{\tX}$ over $Y$.
If one thinks about the structure complexes as $\cO_Y$-algebras, the
existence of such an isomorphism would say that the structure complex
$\cO_{\tX}$ of $\tX$ would be quasi-isomorphic to its cohomology. This is
equivalent to saying that the derived scheme $\tX$ is formal over $Y$
in the sense of~\cite{DGMS}. See~\cite{ACH} for more discussions.

\paragraph
The HKR isomorphism on Hochschild cohomology is obtained by dualizing
the structure complexes of $\bT_X[-1]$ and $LX$ with respect to $X$,
and taking global sections:
\begin{align*}
  \HH^*(X) & = \Hom_{X\times X}(\Delta_*\cO_X,\Delta_*\cO_X) =
        \Hom_X(\Delta^*\Delta_*\cO_X, \cO_X)  \\
       & = \Gamma(X, (\cO_LX)^\chk)  \iso \Gamma(X, \Sym T_X[-1]) =
          \bigoplus H^p(X,\wedge^qT_X) \\
       & = \HT^*(X).
\end{align*}

Duals of functions are distributions.  We will need a relative version
of this concept, made precise in the following definition.

\begin{Definition}
  For a map of spaces $f: Y\ra X$ the space of relative
  distributions is defined by
  \[ \bD(Y/X) = \Hom(f_*\cO_Y,\cO_X). \]
  We will often omit the space $X$ when it is clear from context.
\end{Definition}
\bigskip

\paragraph
For example consider the map $LX\ra X$, where $X$ is a smooth
algebraic variety.  Then the Hochschild cohomology of $X$ is naturally
identified with the space of distributions on $LX$,
\begin{align*}
&  \bD(LX/X)  = \Hom(\Delta^*\cO_\Delta,\cO_X)= \HH^*(X), 
\intertext{while the polyvector field cohomology of $X$ is naturally the
space of distributions on $\bT_X[-1]$, the linearization of $LX$,}
& \bD(\bT_X[-1]) =\Hom(\Sym\Omega_X[1],\cO_X) = \HT^*(X).
\end{align*}
Therefore we should think of polyvector fields as (invariant)
distributions on the Lie algebra $\bT_X[-1]$ and Hochschild cohomology
as (invariant) distributions on the derived group $LX$.  The HKR
isomorphism is then interpreted as the isomorphism on distributions
induced by the exponential map.  The product structures on the two
sides are given by convolution of distributions, where the group
structure on $\bT_X[-1]$ is given by addition in the fibers.

This interpretation was probably the basis for Kontsevich's claim: for
ordinary Lie algebras a theorem of Duflo~\cite{D} asserts that the
rings of invariant distributions on a Lie group and on its Lie
algebra are isomorphic, after a correction to the exponential map by
what is known as the Duflo element.

\paragraph{\bf The Theorem of Kontsevich and Calaque-Van den
  Bergh.}
\label{subsec Kont Thm}
Calaque and Van den Bergh~\cite{CV}, following the claim of
Kontsevich~\cite{K}, proved that the map
\[ \HT^*(X) \stackrel{\Td^{1/2}\lrcorner}{\lra} \HT^*(X)
  \stackrel{\HKR}{\lra} \HH^*(X) \]
is a {\em ring} isomorphism.  Here the analogue of the Duflo element
is the characteristic class $\Td^{1/2}$.

We are interested in studying  generalizations of the above theorem
to the case of global quotient orbifolds.

\paragraph{\bf The orbifold HKR isomorphism.}
Before we begin we need an analogue of the HKR isomorphism for
orbifolds.  Let $X$ be a smooth algebraic variety, let $G$ be a finite
group acting on $X$, and denote by $\sX=[X/G]$ the corresponding
global quotient orbifold.

We have the following diagram
\[ 
\xymatrix{
L\sX\ar[ddr]_p\ar[drr]^q
& & \\ &I\sX\ar[ul]\ar[r]\ar[d] &
\sX\ar[d]^\Delta \\
&\sX\ar[r]^{\Delta} &
\sX\times\sX. }
\] 
Here $L\sX$ denotes the loop space of the stack $\sX$ defined by
analogy with the case of ordinary spaces as the derived
self-intersection
\[  L\sX=\sX\times^R_{\sX\times\sX}\sX. \]
Its underlying underived stack $I\sX$ is the
inertia stack of $\sX$,  
\[ I\sX=\sX\times_{\sX\times\sX}\sX.\]
Unlike the case where $X$ is a smooth space with no group action,
the two maps $I\sX \ra \sX$ are no longer isomorphisms: it is not hard to
see that 
\[ I\sX=\left [(\coprod_{g\in G} X^g)/G \right], \]
where $X^g$ denotes the fixed locus of the action of $g\in G$ on $X$.

We can rewrite $X^g$ as $\Delta\times_{X\times X}\Delta^g$, where
$\Delta=\{(x,x)\}\subset X\times X$ and
$\Delta^g=\{(x,gx)\}\subset X\times X$.  We get an explicit
formula for the derived loop space $L\sX$ if we replace the above intersection
by the corresponding derived intersection:
\[ L\sX = \left [ (\coprod_{g\in G}\tXg)/G\right], \]
where
\[ \tXg=\Delta\times^R_{X\times X}\Delta^g. \]
We will call $\tXg$ the derived fixed locus of $g$.

\paragraph
The orbifold HKR isomorphism expresses the derived loop space
$L\sX$ as the total space of a certain vector bundle over $I\sX$, as
explained by Arinkin, C\u{a}ld\u{a}raru and Hablicsek~\cite{ACH}. 

The derived loop space $L\sX$ decomposes naturally into connected
components, so it is better to look at each component
$\tXg$ of $L\sX$ individually.  The orbifold HKR
isomorphism identifies $\tXg$ with the total space of the
tangent bundle of $X^g$.  More precisely, for each $g\in G$~\cite{ACH}
construct a linearization isomorphism of derived schemes over $X$
\[ \bT_{X^g}[-1] = \bL\tXg \stackrel{\sim}{\lra} \tXg.  \]

In explicit terms this translates into an isomorphism of commutative
$\cO_X$-algebras 
\[  q_*\cO_{\tXg} \stackrel{\sim}{\lra}
  i_{g*}\Sym(\Omega_{X^g}[1]), \]
where $i_g: X^g \hookrightarrow X$ is the inclusion of the fixed locus.
Applying $\Hom(-,\cO_X)$ to this algebra isomorphism we get an
induced isomorphism on distributions
\[ \bD(\bT_{X^g}[-1]/X)\stackrel{\sim}{\lra} \bD(\tXg/X) . \]
We will denote $\bD(\bT_{X^g}[-1]/X)$ by $\HT^*(X; g)$, and $\bD(\tXg/X)$
by $\HH^*(X;g)$.

Grothendieck duality allows us to give an explicit form to the space
$\bD(\bT_{X^g}[-1]/X)$:
\[ \bD(\bT_{X^g}[-1]/X) = \Hom_X(i_{g*} \Sym \Omega_{X^g}[1], \cO_X) =
  \bigoplus_{p+q=*}H^{p-c_g}(X^g,\wedge^qT_{X^g}\otimes\omega_g), \]
where $c_g$ is the codimension of $X_g/X$ and $\omega_g$ is the
dualizing sheaf of the inclusion $X^g\subseteq X$.  Taking
$G$-invariants of the direct sum over $g\in G$ we get the final 
form of the orbifold HKR isomorphism for $\sX$:
\begin{align*}
  \HH^*(\sX) & = \left(\bigoplus_{g\in G}\HH^*(X,g)\right)^G
               =\left (\bigoplus_{g\in G}\bD(\tXg/X)\right)^G  \\
             & \iso \left(\bigoplus_{g\in G} \HT^*(X;g)\right)^G =
               \left (\bigoplus_{g\in G}\bigoplus_{p+q=*} H^{p-c_g}(X^g,\wedge^q
               T_{X^g}\otimes\omega_g) \right)^G.
\end{align*}                                 
We think of the right hand side above as the definition of the space
of polyvector fields on $\sX$,
\[ \HT^*(\sX)  = \left(\bigoplus_{g\in G}  \HT^*(X;g)\right)^G. \]

\paragraph
We close this section by noting that the above HKR isomorphisms can be
assembled to an analogue of the exponential map
\[ \exp:  \bT_{I\sX}[-1] = \bL_{L\sX}\stackrel{\sim}{\lra} \bL\sX \]
from the Lie algebra $\bL_{L\sX}$ to the derived group $L\sX$.

\section{Definition of the product on orbifold polyvector fields}

In this section we define the product on orbifold polyvector fields.
The technical results used in the definition are introduced
in this section, but will be proved in Section 4.

As we have explained previously, the Hochschild cohomology and the
polyvector fields cohomology of a space $X$ can be viewed as the
distributions on the derived loop space (a derived Lie group) and on
its Lie algebra, respectively.  The product structures on these come
from the convolution of distributions.  We begin by recalling the
definition of the convolution product of distributions on (classical)
Lie groups and Lie algebras.

\paragraph{\bf Distributions on Lie groups and Lie algebras.}  Let $G$
be a Lie group with Lie algebra $\g$.  The convolution product of
distributions $\bD(G)$ on $G$ is defined as follows
$$ \xymatrix{ \bD(G)\otimes\bD(G)\ar[r] &
  \bD(G\times G)\ar[r]^{m_*} &\bD(G), }$$ where $m$ is the
multiplication map $G\times G\rightarrow G$, and $m_*$ is the induced
map on distributions.  The Lie algebra $\g$ of $G$ is a vector space.
It is considered as an abelian group under the addition operation of
vectors.  One can define the convolution product for $\bD(\g)$
similarly.

In the derived setting, the convolution product on orbifold Hochschild
cohomology is known as the composition of morphisms in the derived
category.  We hope to define the convolution product on the polyvector
fields.  Therefore, it is important to know how to recover the
convolution product on $\bD(\g)$ with the knowledge of the group $G$
only.  The following is how we do this.

First there is a multiplication map $m:G\times G\rightarrow G$.
Taking derivative of this map, we get the induced map on tangent
spaces $\bL_m: \bL_{G\times G}\rightarrow \bL_G$, where $\bL_G=\g$
is the tangent space of $G$ at origin.  We use the same notation $\bL$
as the notation for the linearization of derived schemes in the
derived setting.

There is a natural isomorphism $\bL_{G\times
G}\cong\bL_G\times\bL_G$.  Under this natural identification, the map
$\bL_m:\bL_G\times\bL_G\rightarrow\bL_G$ is nothing but the addition
law on the vector space $\bL_G$.  We can recover the convolution
product of $\bD(\bL_G)$ now $$\xymatrix{
\bD(\bL_G)\otimes\bD(\bL_G)\ar[r] &
\bD(\bL_G\times\bL_G)\ar[r]^{\cong} &\bD(\bL_{G\times
G})\ar[r]^{{\bL_m}_*} & \bD(\bL_G). } $$

\paragraph{\bf A non-trivial isomorphism in the derived setting.} We
can try to do exactly the same thing in the derived setting.  However,
there is a technical issue.  The natural isomorphism $\bL_{G\times
G}\cong\bL_G\times\bL_G$ is not at all obvious for the derived loop
space.  The analogous statement would be
\[ \bL_{L\sX\times_{\sX}^R  L\sX}\cong\bL_{L\sX}\times^R_{\sX}\bL_{L\sX} \]
for the derived loop space of an orbifold $\sX$.  The left hand side
is obviously linear: it is a total space of a vector bundle over the
inertia stack $I\sX$.  On the other hand, it is not at all obvious
that the right hand side can be linearized.

The following two propositions will be proved in the next sections.

\begin{Proposition}
  \label{Prop formality}
  Let $\sX=[X/G]$ be a global quotient orbifold of a finite group $G$
  acting on a smooth algebraic variety.  If we further assume $G$ is
  abelian, then there is an isomorphism
  \[ \bL_{(L\sX\times_{\sX}^R L\sX)}\cong\bL_{L\sX}\times^R_{\sX}\bL_{L\sX}. \]
\end{Proposition}

The derived loop space $L\sX$ decomposes naturally into connected
components, so we can restate the above proposition on components.

\begin{Proposition}
  \label{Prop formality'}
  In the same setting as Proposition \ref{Prop formality}, there is an
  isomorphism
  \[
    \bL_{\widetilde{X^g}\times^R_X\widetilde{X^h}} \cong
    \bL_{\widetilde{X^g}} \times^R_X\bL_{\widetilde{X^h}} \] 
  for any $g,h\in G$.
\end{Proposition}
\medskip

\paragraph{\bf The definition of the convolution product in the
  derived setting.}
The multiplication map for Lie groups plays an important role in the
case of Lie groups and Lie algebras.  We need to know what the
multiplication map is for the derived loop space
$L\sX=\sX\times^R_{\sX\times\sX}\sX$ of $\sX$.  It is the projection
map $p_1\times p_3$ onto the first and the third factors
\[ \xymatrix{
    L\sX\times_{\sX}L\sX\ar[r]^{m}\ar[d]^{=} & L\sX\ar[d]^{=}.\\
    \sX\times^R_{\sX\times\sX}\sX\times^R_{\sX\times\sX}\sX\ar[r]^{~~~~~~p_1\times
      p_3} & \sX\times^R_{\sX\times\sX}\sX. } \]

We need three lemmas for derived groups which are generalizations of
well-known results from classical Lie group theory.

\begin{Lemma}
  \label{lem induced map on linearizations}
  A map $f: X\longrightarrow Y$ between derived schemes induces a map on
  linearlizations $\bL_f:\bL_{X}\longrightarrow\bL_Y$.
\end{Lemma}

\begin{proof} We have a commutative diagram $$ \xymatrix{ X^0
\ar[d]^{i} \ar[r]^{g} & Y^0 \ar[d]^{j} \\ X \ar[r]^{f} & Y, } $$

where $X^0$ and $Y^0$ are the classical schemes of $X$ and $Y$
respectively.  Then there is a commutative diagram of derived tangent
complexes $$ \xymatrix{ T_{X^0} \ar[d] \ar[r] & g^*T_{Y^0} \ar[d] \\
i^*{T}_{X} \ar[r] & i^*f^*{T}_{Y}=g^*j^*{T}_{Y}.  } $$

Passing to the quotient, we get an induced
map $$N_{X^0/X}=i^*{T}_{X}/T_{X^0}\rightarrow
g^*(j^*{T}_{Y})/g^*T_{Y^0}=g^*(j^*{T}_{Y}/T_{Y^0})=g^*N_{Y^0/Y}.$$

The map above is equivalent to a map $\bN_{X^0/X}\rightarrow
\bN_{Y^0/Y}\times_{Y^0}X^{0}$ in terms of total spaces. \end{proof}

Applying the above lemma to the multiplication map of derived loop
space yields an induced map
$\bL_m:\bL_{L\sX\times^R_{\sX}L\sX}\rightarrow\bL_{L\sX}.$

\begin{Lemma}
  \label{lem push-forward of dist}
  Suppose there is a commutative diagram
  $$\xymatrix{ X \ar[rr]^f\ar[dr]_i & & Y\ar[dl]^j \\ & S & }$$ of
  (derived) schemes.  Then there is a pushforward map for relative
  distributions, i.e., there is a natural induced map
  $f_*:\bD(X/S)\rightarrow\bD(Y/S).$
\end{Lemma}

\begin{proof}
  We have $\bD(X/S)=\Hom(i_*\cO_X,\cO_S)$.  Applying $j_*$ to the map
  $\cO_Y\rightarrow f_*\cO_X$, we get a map
  $j_*\cO_Y\rightarrow j_*f_*\cO_X=i_*\cO_X$.  Composing it with
  $i_*\cO_X\rightarrow \cO_S$, we get the desired pushforward map.
\end{proof}

\begin{Lemma}
  \label{lem tensor prod of dist}
  Suppose there is a commutative diagram
  $$ \xymatrix{ W=X\times^R_SY\ar[r]\ar[d]\ar[dr]^\pi & Y\ar[d]_j \\
    X\ar[r]^i & S} $$ of (derived) schemes.  Then there is a natural map
  $\bD(X/S)\otimes\bD(Y/S)\rightarrow\bD(W/S)$.
\end{Lemma}

\begin{Proof}
  We have
  \begin{align*}
    \Hom(i_*\cO_X,\cO_S)\otimes\Hom(j_*\cO_Y,\cO_S) & \ra
    \Hom(i_*\cO_X\otimes_{\cO_S}   j_*\cO_Y,\cO_S) \\
                                                    &
                                                      =\Hom(\pi_*\cO_W,\cO_S). 
  \end{align*}
  \vspace*{-4em}
  
\end{Proof}

With the three lemmas above we are able to define our desired
product.

\begin{Definition}
\label{def product}
Under the assumptions in Proposition~\ref{Prop formality} we define
the following binary operation on $\bD(\bL_{L\sX}/\sX)$, which is our
proposed definition for a product on orbifold polyvector fields:.
\begin{align*}
  \bD(\bL_{L\sX}/\sX)\otimes\bD(\bL_{L\sX}/\sX) & \lra
      \bD(\bL_{L\sX}\times^R_{\sX}\bL_{L\sX}/\sX) \stackrel{\sim}{\lra}
      \bD(\bL_{(L\sX\times_{\sX} L\sX)}/\sX) \\
      & \stackrel{\bL_{m_*}}{\lra}\bD(\bL_{\sX}/\sX), 
\end{align*}
where the first arrow is due to Lemma \ref{lem tensor prod of dist},
the second arrow is the non-trivial isomorphism in Proposition
\ref{Prop formality}, and the last map is due to Lemmas \ref{lem
  induced map on linearizations} and \ref{lem push-forward of dist}.
\end{Definition}
\medskip

Looking at each connected component of $L\sX$ individually, the
definition gives a map for every $g, h\in G$
\[
  \bD(\bL_\tXg/X)\otimes\bD(\bL_\tXh/ X) \ra
    \bD(\bL_\tXg\times^R_{X}\bL_\tXh/
    X)\stackrel{\sim}{\lra} \bD(\bL_{(\tXg\times_{X}
      \tXh)}/ X)\stackrel{\bL_{m_*}}{\lra}
    \bD(\bL_{\tXgh}/ X). 
\]

\section{The formality of double fixed loci}

We begin by studying the cohomology sheaves of the structure complex
of $\bL_{\tXg}\times^R_X\bL_{\tXh}$.  Then we compute the
linearization $\bL_{\tXg\times^R_X \tXh}$ explicitly.  At last, we
prove Propositions~\ref{Prop formality} and~\ref{Prop formality'}, in
other words we construct a formality isomorphism
\[ \bL_{\tXg}\times^R_X\bL_{\tXh}\cong \bL_{\tXg\times^R_X \tXh}. \]
The construction will be indirect: we will use known results to show
that both sides are isomorphic to
\[ \bT_{X^g}[-1]|_{X^{g,h}} \oplus \bT_{X^h}[-1]|_{X^{g,h}} \oplus \bE[-1], \]
where $\bE$ is the total space of the excess intersection bundle for
the intersection of $X^g$ and $X^h$ in $X$.

Throughout this section $X$ will be a smooth variety, and $g,h$ will denote
commuting elements of a group $G$ which acts on $X$.

\paragraph
Before we begin we note that the derived fixed locus
$\tXg=\Delta\times^R_{X\times X}\Delta^g$ is not directly a scheme
over $X$.  It is more naturally viewed as a scheme over $\Delta$ or
$\Delta^g$, both of which are isomorphic (in different ways) to $X$.
We will use the latter when computing the fiber product
$\tXg\times_X^R \tXh$.

Similarly, the derived fixed locus
$\tXh=\Delta\times^R_{X\times X}\Delta^h$ is naturally isomorphic to
$\Delta^g\times^R_{X\times X}\Delta^{gh}$, which is also a scheme over
$\Delta^g$.  Therefore, while the notation $\tXg\times_X^R\tXh$ is
imprecise, what we will really mean by it is
\[ (\Delta\times^R_{X\times
    X}\Delta^g)\times_{\Delta^g}^R(\Delta^g\times_{X\times
    X}^R\Delta^{gh})=\Delta\times^R_{X\times X}\Delta^g\times^R_{X\times
    X}\Delta^{gh}.\]
We think of this as the derived fixed locus of $g$ and $h$, and denote
it by $\tXgch$.

\paragraph{\bf The cohomology sheaves of the structure complex of $\bL_\tXg \times_X^R \bL_\tXh$.}
It is difficult to compute $\cO_{\bL_{\tXg}\times^R_X\bL_{\tXh}}$
directly, but we can compute its cohomology sheaves more easily, and
we begin with this computation.

We hope to compute the cohomology sheaves of
\[ \cO_{\bL_\tXg \times_X^R \bL_\tXh} = \Sym(\Omega_{X^g}[1])\otimes^{L}_{\cO_X}
  \Sym(\Omega_{X^h}[1]). \]
The calculation becomes straightforward using the following lemma.

\begin{Lemma}
  Suppose $i$, $j$ are closed embeddings of classical schemes, and
  $\cE$, $\cF$ are vector bundles on $X$ and on $Y$, respectively.
  Denote by $W$ the fiber product of classical schemes below
  \[
    \xymatrix{ X\times_S Y=W \ar[d]_{k} \ar[r]^{~~~~~~l} & Y \ar[d]^{j} \\
      X \ar[r]^{i} & S.  }
  \]
  Then
  \[ \sH_q(i_*\cE\otimes^L_{\cO_S}j_*\cF)=j_*l_*(\cE|_{W}\otimes\cF|_{W}\otimes
    E^{\vee}), \]
  where $E$ is the excess intersection bundle, 
  \[ E=\frac{T_S|_W}{T_X|_W+T_Y|_W}. \]  
\end{Lemma}

\begin{Proof}
  \cite[Proposition A.6]{CKS}.
\end{Proof}

\paragraph
The lemma above shows that
\[ \sH^n(\cO_{\bL_{\tXg}\times^R_X\bL_{\tXh}}) =
  \bigoplus_{p+q+i=n}(\wedge^p\Omega_{X^g}|_{X^{g,h}} \otimes
  \wedge^q\Omega_{X^h}|_{X^{g,h}} \otimes \wedge^i E^{\vee}). \]
In other words, if we knew that $\bL_\tXg \times_X^R \bL_\tXh$ is formal over its underlying
classical scheme $X^{g,h}$, the above calculation would imply that
\[ \bL_\tXg \times_X^R \bL_\tXh \iso \bT_{X^g}[-1]|_{X^{g,h}} \oplus
  \bT_{X^h}[-1]|_{X^{g,h}} \oplus \bE[-1].\]
We will prove the formality statement in~(\ref{subsec:formal}).

\paragraph{\bf The structure complex of $\bL_{\tXgch}$.}
\label{4.3}
We will now argue that the linearization $\bL_{\tXgch}$
is precisely the space that appears on the right hand side of the
equality above,
\[ \bL_{\tXgch} = \bT_{X^g}[-1]|_{X^{g,h}} \oplus
  \bT_{X^h}[-1]|_{X^{g,h}} \oplus \bE[-1]. \]

By definition,
$\bL_{\tXgch}=\bN_{X^{g,h}/\tXgch}$.
To compute the normal bundle, we need to know what the derived tangent
complex of $\tXgch$ is.

\paragraph{\bf The derived tangent complex.}
\label{derived tangent complex}
The standard reference for the derived tangent complex is~\cite{I}.
Suppose $X$ and $Y$ are closed subschemes of a scheme $S$.  Let
$\tW=X\times^R_SY$ be the derived intersection and $W=X\times^R_SY$ be
the classical intersection.  There is an exact sequence
$$0\rightarrow T_X|_W\cap T_Y|_W=T_W\rightarrow T_X|_W\oplus T_Y|_W\rightarrow T_S|_W\rightarrow E\rightarrow0.$$

The complex
\[ T_X|_W\oplus T_Y|_W\rightarrow T_S|_W=\Cone(T_X|_W\oplus T_Y|_W
  \rightarrow T_S|_W)[-1] \]
is the restriction to $W$ of the derived tangent complex ${T}_{\tW}|_W$ of 
$\tW$.  Since only its $\sH^0$ and $\sH^1$ sheaves are non-zero, the
information contained in it is equivalent to the data of 
the triple $(\sH^0, \sH^1, \eta)$, where $\sH^0({T}_{\tW}|_W)=T_W$,
$\sH^1({T}_{\tW}|_W)=N_{W/\tW}=E$, and the class $\eta$ is an element
in $\Ext^2_S(E,T_W)$.

For example, if we consider the situation where $S=X\times X$,
$X=\Delta$, $Y=\Delta^g$, so that $W=X^g$, we have
$\sH^0({T}_{\tXg})=T_{X^g}$ and $\sH^1({T}_{\tXg})=E$.  Moreover, the
excess bundle in this case equals the coinvariant bundle
$(T_X|_{X^g})_g$, which in characteristic zero is canonically
isomorphic to the invariant bundle $T_{X^g}$.

The linearization $\bL_{\tW}$ is by definition the total space of the
normal bundle $N_{W/\tW}$, the cone of the map $T_W \ra T_\tW|_W$.
Since $\sH^0(T_\tW|_W) = T_W$, it follows that the normal bundle
$N_{W/\tW}$ is the first cohomology $\sH^1({T}_{\tW}|_W)[-1]$ of
$T_{\tW}|_W$.  In the example considered above this shows that
$\bL_{\tXg} = \bT_{X^g}[-1]$.

\paragraph
The above discussion also works for derived schemes.  If we replace
$X$, $Y$, and $S$ by derived schemes in the commutative diagram at the
beginning of~(\ref{derived tangent complex}), we have the same formula
$$T_{\tW}|_W=\Cone({T}_{X}|_W\oplus {T}_{Y}|_W\longrightarrow {T}_S|_W)[-1],$$
where $T_X$, $T_Y$, and $T_S$ are the derived tangent complexes of
$X$, $Y$, and $S$. The scheme $W$ is the underlying classical scheme
of $\tW=X\times^R_SY$.  Since all the complexes are
restricted to $W$, we will omit the restrictions from $X$, $Y$, $S$,
and $\tW$ to $W$ for simplicity from now on.

It helps us to compute the derived tangent complex
$T_{\tXgch}$ if we set $S=X$, $X=\tXg$, $Y=\tXh$, and
$\tW=\tXgch$ respectively.

\begin{Lemma}
The derived tangent complex of $\tXgch$ is quasi-isomorphic to
\[ T_{\Delta}\oplus T_{\Delta^g}\oplus T_{\Delta^h}\rightarrow
T_{X\times X}\oplus T_{X\times X}\oplus T_{X\times X}\rightarrow
T_{X\times X}, \]
where the maps are of the form
\begin{align*}
  (v_1,v_2,v_3) & \rightarrow(v_1-v_2,v_2-v_3,v_3-v_1), 
\intertext{and}
                  (a,b,c) &\rightarrow a+b+c.
\end{align*}
\end{Lemma}

We can compute the cohomology of the derived tangent complex of
$\tXgch$ using the above lemma.  We have
$\sH^0({T}_{\tXgch})=T_{X^{g,h}}$. To compute the first cohomology it
suffices to compute the cokernel of the map below
\[ V\oplus V\oplus V\rightarrow V\oplus V\oplus V\oplus V, \]
where $(V=T_{X}\cong T_{\Delta}\cong T_{\Delta^g}\cong T_{\Delta^h})$ and
the maps are $(v,v',v'')\rightarrow(v-v',v-gv',v-v'',v-hv'')$.  This
is done in the lemma below.

\begin{Lemma}
  Suppose $V$ is a finite dimensional representation of a finite group
  $G$ over a field of characteristic $0$.  Let $g$ and $h$ be two elements of $G$.
  Then the quotient of $V\oplus V\oplus V\oplus V$ by the
  relations $(v,v,v,v)$, $(v,gv,0,0)$, and $(0,0,v,hv)$ is isomorphic
  to
  \[ V_g\oplus V_h\oplus \frac{V}{V^g+V^h}. \]
\end{Lemma}

\begin{proof}
  Let $L$ be the linear subspace $(v,v,v,v)$, and note that
  $L=H_1\cap H_2\cap H_3$, where $H_{1}$ is defined by $v_1=v_2$,
  $H_{2}$ is defined by $v_1=v_3$, and $H_3$ is defined by $v_3=v_4$.
  Then we have an isomorphism
\[ \frac{V\oplus V\oplus V\oplus V}{L}\cong \frac{V^{\oplus
      4}}{H_1}\oplus\frac{V^{\oplus4}}{H_2}\oplus\frac{V^{\oplus4}}{H_3}\cong
  V\oplus V\oplus V.\]

Under this identification, the second and third relations become
$(v-gv,v,0)$ and $(0,-v,v-hv)$.  There is a natural projection to the
first and third components,
\[\frac{V\oplus V\oplus V}{(v-gv,v,0),(0,-v',v'-hv')}\rightarrow
  \frac{V\oplus V}{(v-gv,0),(0,v'-hv')}=V_g\oplus V_h.\]
It is easy to show that the kernel is $\frac{V}{V^g+V^h}$.  So we get a
short exact sequence
\[ 0\rightarrow\frac{V}{V^g+V^h}\rightarrow\frac{V\oplus V\oplus
  V}{(v-gv,v,0),(0,-v',v'-hv')}\rightarrow V_g\oplus V_h\rightarrow 0. \]

By the averaging map
$v\rightarrow \frac{1}{\ord(g)}\sum_{i=1}^{\ord(g)}g^i\cdot v$, the map
$V\rightarrow V_g$ splits in characteristic $0$.  We can use the
averaging map of $g$ and $h$ to get a canonical splitting of the short
exact sequence above.
\end{proof}

The discussion above shows that the first cohomology of the tangent
complex of $\tXgch$ is $E\oplus T_{X^g}\oplus
T_{X^h}$, where $E=\frac{T_{X}}{T_{X^g}+T_{X^h}}$.  As a consequence
we have an isomorphism
\[ \bL_{\tXg\times_X^R \tXh} = \bL_\tXgch\iso \bT_{X^g}|_{X^{g,h}}[-1] \oplus
  \bT_{X^h}|_{X^{g,h}}[-1] \oplus \bE[-1]. \]

\paragraph{\bf Formality of $\bL_{\tXg}\times^R_X\bL_{\tXh}$.}
\label{subsec:formal}
The linearizations $\bL_{\tXg}$ and $\bL_{\tXh}$ are by definition the
total spaces of vector bundles over $X^g$ and $X^h$ respectively, so we
study the formality of $X^g\times^R_X X^h$ first.  The key tools are
\cite[Theorem 1.8 and Lemma 4.3]{ACH}.

\begin{Proof}[Proof of Proposition~\ref{Prop formality'}.]  The
inclusion $X^g\rightarrow X$ splits to first order.  By \cite[Lemma
4.3]{ACH}, the derived scheme $X^g\times^R_X X^h$ is formal over
$X^g\times X^h$ if and only if the short exact sequence on
$X^{g,h}=X^g\times_X X^h$
$$0\rightarrow\frac{T_{X^h}}{T_{X^g}\cap T_{X^h}}\rightarrow    \frac{T_X}{T_{X^g}}\rightarrow E=\frac{T_X}{T_{X^g}+T_{X^h}}\rightarrow0$$
splits.

Define a map
\begin{align*}
  \frac{T_X}{T_{X^g}} & \rightarrow\frac{T_{X^h}}{T_{X^g}\cap T_{X^h}}
                        \intertext{by the formula}
                        v & \mapsto \frac{1}{\ord(h)}\sum h^i\cdot v.
\end{align*}
The map is well-defined because $g$ and $h$ commute under our initial
assumptions (or the ones in Proposition~\ref{Prop formality'}).  It
splits the short exact sequence above.  This shows $X^g\times^R_X X^h$
is formal over $X^g\times X^h$.

Consider the following commutative diagram
\[
\xymatrix{ \widehat{X^{g,h}}=X^g\times^R_X X^h \ar[drr]\ar[ddr] & & \\
& X^{g,h}\ar[ul] \ar[d]^q \ar[r]_{p} & X^g \ar[d]_i \\ & X^h \ar[r]^j
& X }
\]
By \cite[Theorem 1.8]{ACH} we know that the dg functor $j^*i_*(-)$ is
isomorphic to $q_*\left (p^*(-)\otimes\Sym (E^\vee[1])\right )$.

The structure complex of $\bL_{\tXg}\times^R_X\bL_{\tXh}$ is
\[ i_*\Sym(\Omega_{X^g}[1])\otimes^L_{\cO_X}j_*\Sym(\Omega_{X^h}[1]) =
  j_*\left (j^*i_*\Sym(\Omega_{X^g}[1])\otimes \Sym(\Omega_{X^h}[1])\right ). \]

Using the isomorphism of the two dg functors above, we see that
$j^*i_*(\Sym(\Omega_{X^g}[1]))\cong
q_*(p^*(\Sym(\Omega_{X^g}[1]))\otimes\Sym(E^\vee[1]))$. 
As a consequence
\begin{align*}
  i_*\Sym(\Omega_{X^g}[1]) & \otimes^L_{\cO_X}j_*\Sym(\Omega_{X^h}[1])
  = \\
   & =
     j_*q_*\left (\Sym(\Omega_{X^g}|_{X^{g,h}}[1])\otimes\Sym(E^\vee[1])
     \otimes\Sym(\Omega_{X^h}|_{X^{g,h}}[1]) \right )  \\
  & = j_*q_*\Sym\left ((\Omega_{X^g}|_{X^{g,h}}\oplus\Omega_{X^h}|_{X^{g,h}}
    \oplus E^\vee )[1] \right ).
\end{align*}
Therefore $\bL_{\tXg}\times^R_X\bL_{\tXh}$ is formal over $X$ (and it is
isomorphic to $\bL_{\tXgch}$).
\end{Proof}

\section{Associativity of the product}

In this section we explain the strategy for proving Theorem A.  In
other words we want to show that, under the assumption that certain
Bass-Quillen cohomology classes vanish, the product defined in Section
3 is associative. The proof is reduced to Propositions~\ref{prop assoc
  1} and~\ref{prop assoc 2}, which will be proved in Section 6.

\paragraph{\bf Formality of triple intersections.}
To prove the associativity it is natural to study the triple
intersection
$\widetilde{X^{g,h,k}}=\tXg\times^R_X\times\tXh\times^R_X\tXk$
for $g$, $h$, and $k\in G$.  More precisely we define
\begin{align*}
  \widetilde{X^{g,h,k}} & = (\Delta\times^R_{X\times X}\Delta^g)
  \times_{\Delta^g} (\Delta^g\times^R_{X\times X}\Delta^{gh})
  \times_{\Delta^{gh}}(\Delta^{gh}\times^R_{X\times X}\Delta^{ghk}) \\
& =\Delta\times^R_{X\times X}\Delta^g \times^R_{X\times X}
                         \Delta^{gh}\times^R_{X\times X} \Delta^{ghk},
\end{align*}
as explained in~(\ref{4.3}).  Under the assumption that $G$ is abelian
it is not hard to see that $\widetilde{X^{g,h,k}}$ is formal over $X$,
and it is isomorphic to $\bL_{\widetilde{X^{g,h,k}}}$.  The proof is
essentially the same as the one in Section 4.

\paragraph
\label{big diagram}
The diagram
\[
\xymatrix{ \widetilde{X^{g,h,k}}\ar[r]\ar[d] &
\widetilde{X^{g,hk}}\ar[d]\\ \widetilde{X^{gh,k}}\ar[r]&
\widetilde{X^{ghk}} }
\] 
is commutative because it is the associativity of the group law of the
loop space $L[X/G]$.  Taking distributions over on the corresponding
linearizations, we get the following commutative diagram
\[ 
\xymatrix{
\bD(\bL_{\widetilde{X^{g,h,k}}})\ar[r]\ar[d] &
\bD(\bL_{\widetilde{X^{g,hk}}})\ar[d]\\
\bD(\bL_{\widetilde{X^{gh,k}}})\ar[r]&
\bD(\bL_{\widetilde{X^{ghk}}}).  }
\]
For simplicity we have denoted the relative distributions with respect to
$X$ as $\bD(-)$ instead of $\bD(-/X)$.

\paragraph{\bf What we need to prove.} 
Consider the following diagram:
\[ \mbox{\small
\xymatrix{
\bD(\bL_{\tXg})\otimes\bD(\bL_{\tXh})\otimes\bD(\bL_{\tXk})\ar[ddd]_{m\otimes
\id}\ar[rrr]^{\id\otimes m}\ar[dr] & &
&\bD(\bL_{\widetilde{X^{g}}})\otimes\bD(\bL_{\widetilde{X^{hk}}})\ar[dl]\ar[ddd]^{m}
\\ & \bD(\bL_{\widetilde{X^{g,h,k}}})\ar[r]\ar[d] &
\bD(\bL_{\widetilde{X^{g,hk}}})\ar[d] \\ &
\bD(\bL_{\widetilde{X^{gh,k}}})\ar[r]&
\bD(\bL_{\widetilde{X^{ghk}}}) & \\
\bD(\bL_{\widetilde{X^{gh}}})\otimes\bD(\bL_{\widetilde{X^{k}}})\ar[rrr]_{m}\ar[ur]
& & & \bD(\bL_{\widetilde{X^{ghk}}})\ar[ul]^{=}, }} 
\] 
where $m$ is the product on orbifold polyvector fields in
Definition~\ref{def product}.  Associativity of $m$ is equivalent to
commutativity of the outer part of the diagram.

The middle square is commutative by the discussion in the previous
paragraph.  The squares on the bottom and right are commutative
because they are the definitions of our product.

We need to examine the ones on the top and left.  The left one is the
diagram
\[ 
  \xymatrix{\bD(\bL_{\tXg})\otimes\bD(\bL_{\tXh})\otimes
   \bD(\bL_{\tXk})\ar[r]\ar[d]_{m\otimes \id} & \bD(\bL_{\widetilde{X^{g,h,k}}})\ar[d] \\
\bD(\bL_{\widetilde{X^{gh}}})\otimes\bD(\bL_{\widetilde{X^{k}}})\ar[r]
& \bD(\bL_{\widetilde{X^{gh,k}}}).  }
\]

Expand the diagram in detail
\[\xymatrixcolsep{0.15in}
\xymatrix{\bD(\bL_{\tXg})\otimes\bD(\bL_{\tXh})\otimes\bD(\bL_{\tXk})\ar[r]\ar[d]_{m\otimes
\id}
&\bD(\bL_{\widetilde{X^{g,h}}})\otimes\bD(\bL_{\widetilde{X^{k}}})\ar[d]\ar[r]
&\bD(\bL_{\widetilde{X^{g,h}}}\times^R_X\bL_{\widetilde{X^{k}}})\ar[r]^{~~~~~\sim}\ar[d]
& \bD(\bL_{\widetilde{X^{g,h,k}}})\ar[d] \\
\bD(\bL_{\widetilde{X^{gh}}})\otimes\bD(\bL_{\widetilde{X^{k}}})\ar[r]^{=}
&
\bD(\bL_{\widetilde{X^{gh}}})\otimes\bD(\bL_{\widetilde{X^{k}}})\ar[r]
&
\bD(\bL_{\widetilde{X^{gh}}}\times^R_X\bL_{\widetilde{X^{k}}})\ar[r]^{~~~~~\sim}&
\bD(\bL_{\widetilde{X^{gh,k}}}).  }
\]

Clearly, the left and middle squares of the diagram above are
commutative.  We have to show commutativity of the square on the
right.  Note that the maps on distributions are induced from maps on
spaces, so we only need to show the commutativity of the diagram below
\[
\xymatrix{
\bL_{\widetilde{X^{g,h}}}\times_X^R\bL_{\widetilde{X^{k}}}\ar[r]^{~~~~~\sim}\ar[d]&
\bL_{\widetilde{X^{g,h,k}}}\ar[d]
\\
\bL_{\widetilde{X^{gh}}}\times_X^R\bL_{\widetilde{X^{k}}}\ar[r]^{~~~~~\sim} 
&\bL_{\widetilde{X^{gh,k}}}. }
\]

Similarly, the commutativity of the top square in the big diagram in~(\ref{big
diagram}) reduces to the commutativity of the diagram below
\[
\xymatrix{
\bL_{\widetilde{X^{g}}}\times_X^R\bL_{\widetilde{X^{h,k}}}\ar[r]^{~~~~~\sim}\ar[d]&
\bL_{\widetilde{X^{g,h,k}}}\ar[d]
\\
\bL_{\widetilde{X^{g}}}\times_X^R\bL_{\widetilde{X^{hk}}}\ar[r]^{~~~~~\sim}
&\bL_{\widetilde{X^{g,hk}}}. }
\]
We will only analyze the former diagram; the proof of the
commutativity of the latter is entirely similar. 

\paragraph
There is, however, one more compatibility that needs to be discussed.
Even though we wrote the top left diagonal map in the big
diagram~(\ref{big diagram}) as a single map, it is in fact clear from
the discussion above that there are two maps here,
\[ \bD(\bL_{\tXg})\otimes\bD(\bL_{\tXh})\otimes\bD(\bL_{\tXk})
  \stackrel{\ra}{\ra}\bD(\bL_{\widetilde{X^{g,h,k}}}). \]
One is the one that appears in the left square in the big diagram
in~(\ref{big diagram}), and the other one is the one that is in the
top square of the big diagram in~(\ref{big diagram}).  We need to
prove that these two maps are the same.  This question is easily
reduced to the following problem.

As mentioned in the previous section, the linearizations $\bL_{\tXg}$,
$\bL_{\tXh}$, and $\bL_{\tXk}$ are vector bundles over the underlying
schemes $X^g$, $X^h$, and $X^k$.  To prove the formality of derived
intersections of linearizations, it suffices to prove the formality of
the underlying schemes.  There are two ways to define the isomorphism
$\bL_{\tXg}\times^R_X\bL_{\tXh}\times^R_X\bL_{\tXk}\cong
\bL_{\widetilde{X^{g,h,k}}}$.  One uses the fact that
$\widehat{X^{g,h}}=X^g\times^R_X X^h$ is formal and
$\widehat{X^{(g,h),k}}={X^{g,h}}\times^R_X X^k$ is formal.  The other
uses the fact that $\widehat{X^{h,k}}=X^h\times^R_X X^k$ is formal
and $\widehat{X^{g,(h,k)}}=X^g\times^R_X{X^{h,k}}$ is formal.
Therefore, we need to prove that the two isomorphisms agree, i.e.,
that the diagram below is commutative
\[\mbox{\scriptsize
\xymatrixcolsep{0.1in} \xymatrix{ (X^g\times^R_X X^h)\times^R_X X^k
\ar[d]^{\id}\ar[r]^{\sim~~~~~~~~~~~~~~~~~~} &
\bL_{\widehat{X^{g,h}}}\times^R_X X^k =
\bL_{\widehat{X^{g,h}}}\times^R_{X^{g,h}}(X^{g,h}\times^R_X
X^k)\ar[r]^{~~~~~~~~~~~~~~\sim}&\bL_{\widehat{X^{g,h}}}\times^R_{X^{g,h}}\bL_{\widehat{X^{(g,h),k}}}
\ar[r] & \bL_{\widehat{X^{g,h,k}}}\ar[d]^{\id}\\ X^g\times^R_X
(X^h\times^R_X X^k) \ar[r]^{\sim~~~~~~~~~~~~~~~~~~} & X^{g}\times^R_X
\bL_{\widehat{X^{h,k}}}=(X^g\times^R_X
X^{h,k})\times^R_{X^{h,k}}\bL_{\widehat{X^{h,k}}}\ar[r]^{~~~~~~~~~~~~~~\sim}
&
\bL_{\widehat{X^{g,(h,k)}}}\times^R_{X^{h,k}}\bL_{\widehat{X^{h,k}}}\ar[r]
& \bL_{\widehat{X^{g,h,k}}}.  } }
\]

\paragraph
\label{def:bq}
Unfortunately we can not prove the commutativity of the diagrams above
without further assumptions.  There is a cohomology class which plays
an important role in what follows.  It appears in a more general
setting, which we review now.

Let $X\hookrightarrow Y\hookrightarrow S$ be a sequence of closed
embedding of smooth schemes, and assume that there is a fixed first
order splitting of the map $X\hookrightarrow Y$, i.e., we have fixed a
map $X^{(1)} \ra X$ which splits the inclusion $X \ra X^{(1)}$ of $X$
into its first order neighborhood $X^{(1)}$ in $Y$.
  
The class we need is the Bass-Quillen class associated to the
restriction $N_{Y/S}|_{X^{(1)}}$ of the normal bundle $N_{Y/S}$ to the
first order neighborhood $X^{(1)}$.  It was introduced by the second
author in~\cite{H}.  In this paper, we call this class the
Bass-Quillen class associated to the sequence of embeddings
$X\hookrightarrow Y\hookrightarrow S$.

The following two statements, which will be proven in the next
section, imply the commutativity of the diagrams above, under the
assumption that the Bass-Quillen classes associated to
$X^{g,h}\hookrightarrow X^{gh}\hookrightarrow X$ and
$X^{g,h}\hookrightarrow X^{g}\hookrightarrow X$ vanish for all
$g,h\in G$.  This will complete the proof of Theorem A.

\begin{Proposition}
  \label{prop assoc 1}
  Under the assumptions of Theorem A, assume that the Bass-Quillen
  class associated to $X^{g,h}\hookrightarrow X^{gh}\hookrightarrow X$
  vanishes.  Then the diagram
  \[
    \xymatrix{
      \bL_{\widetilde{X^{g,h}}}\times_X^R\bL_{\widetilde{X^{k}}}\ar[r]^{~~~~~\sim}\ar[d]&
      \bL_{\widetilde{X^{g,h,k}}}\ar[d] \\
      \bL_{\widetilde{X^{gh}}}\times_X^R\bL_{\widetilde{X^{k}}}\ar[r]^{~~~~~\sim} &
      \bL_{\widetilde{X^{gh,k}}}. }
  \]
  is commutative.
\end{Proposition}

\begin{Proposition}
  \label{prop assoc 2}
  Under the assumptions of Theorem A, assume that the Bass-Quillen
  class associated to $X^{g,h}\hookrightarrow X^{g}\hookrightarrow X$
  and $X^{g,h}\hookrightarrow X^{h}\hookrightarrow X$ vanish.  Then
  the diagram
  \[\xymatrix{
      (X^g\times^R_X X^h)\times^R_X X^k \ar[d]^{\id}\ar[r]^{~~~~~~\sim} &
      \bL_{\widehat{X^{g,h}}}\times^R_X X^k \ar[r]^{~~~~\sim} &
      \bL_{\widehat{X^{g,h,k}}}\ar[d]^{\id}\\ X^g\times^R_X
      (X^h\times^R_X X^k) \ar[r]^{~~~~~~\sim} & X^{g}\times^R_X
      \bL_{\widehat{X^{h,k}}}\ar[r]^{~~~~\sim} &
      \bL_{\widehat{X^{g,h,k}}}.  }
  \]
  is commutative.
\end{Proposition}

\paragraph{\bf Examples.}
\begin{itemize}
\item If $X$ is affine, then all the Bass-Quillen classes above are
zero.
  
\item Consider the $G=\Z/2\Z$ action on an abelian
variety $X$.  We have either $X^{g,h}=X^{gh}$ or $X^{gh}=X$ in this
case, so it is easy to show that all the Bass-Quillen classes are
zero.
\end{itemize}
Therefore the product on $\HT^*(X;G)$ defined in Section 3 is
associative in the cases above.

\section{Consequences of vanishing of Bass-Quillen classes}

We prove Propositions \ref{prop assoc 1} and \ref{prop assoc 2} in
this section.

\paragraph
\label{6.1}
All the linearizations are total spaces of vector bundles over the
underlying schemes, so we can reduce the result of Proposition
\ref{prop assoc 1} to the commutativity of the following two
formality isomorphisms
\[
\xymatrix{ X^{g,h}\ar[d] & X^{g,h}\times_X^R
X^k\ar[r]^{\sim~~~~~~}\ar[d]\ar[l] &
\bE_{(g,h),k}[-1]=\bL_{\widehat{X^{(g,h),k}}}\ar[d] \\
X^{gh}\ar[d]& X^{gh}\times_X^R X^k\ar[d]\ar[r]^{\sim~~~~~~}\ar[l]&
\bE_{gh,k}[-1]=\bL_{\widehat{X^{gh,k}}}\\ X & X^k\ar[l]
& }
\]
where $E_{(g,h),k}=\frac{T_X}{T_{X^{g,h}}+T_{X^k}}$ and
$E_{gh,k}=\frac{T_X}{T_{X^{gh}}+T_{X^k}}$ are excess bundles supported
on $X^{g,h,k}$ and $X^{gh,k}$ respectively.

To check the commutativity, we need to look at how the isomorphism is
defined in \cite{ACH}.  For simplicity, denote $T_X$ by $V$. The two
isomorphisms are defined based on two splittings of the two short
exact sequences
\[
\xymatrix{ 0\ar[r] & \displaystyle{\frac{V^k}{V^{g,h}\cap
V^k}}\ar[r]\ar[d]
&\displaystyle{\frac{V}{V^{g,h}}=N_{X^{g,h}/X}}\ar@{-->}@/_1.6pc/[l]\ar[r]\ar[d]
&\displaystyle{\frac{V}{V^{g,h}+V^k}=E_{(g,h),k}}\ar[r]\ar[d] & 0 \\
0\ar[r] & \displaystyle{\frac{V^k}{V^{gh}\cap V^k}}\ar[r]
&\displaystyle{\frac{V}{V^{gh}}=N_{X^{gh}/X}}\ar@{-->}@/_1.6pc/[l]\ar[r]
&\displaystyle{\frac{V}{V^{gh}+V^k}=E_{gh,k}}\ar[r] & 0 .  }
\]
The two splittings are compatible in the sense that the diagram above
commutes because the two splittings are the averaging map by the
element $k\in G$
\[v\rightarrow\frac{1}{\ord(k)}\sum_{i=1}^{\ord(k)}k^i\cdot v.\]

Proposition \ref{prop assoc 1} is a consequence of the more general
result Proposition \ref{prop assoc 1'} below, by replacing $X$, $Y$,
$Z$, and $S$ in by $X^{g,h}$, $X^{gh}$, $X^k$, and $X$.  Note that all
the assumptions in Proposition \ref{prop assoc 1'} except for the last
one hold trivially for $X^{g,h}$, $X^{gh}$, $X^k$, and $X$.

\begin{Proposition}
  \label{prop assoc 1'}
  Consider a sequence of closed embeddings
  $X\hookrightarrow Y\hookrightarrow S$, and a separate closed
  embedding $Z\hookrightarrow S$.
  
  Assume that all the closed embeddings split to first order in the
  sense of~\cite{ACH}, and that the first order splittings of
  $X\hookrightarrow Y\hookrightarrow S$ are compatible in the sense
  of~\cite{H}.  We further assume that the Bass-Quillen class
  associated to $X\hookrightarrow Y\hookrightarrow S$ is zero. Then the
  diagram
  \[
    \xymatrix{ X\times^R_S Z \ar[r]^{\sim~~~~~~~~~~~~~~}\ar[d] &
      \bE_{X,Z}[-1]=\bE_W[-1]=\bL_{X\times^R_S
        Z}\ar[d]\\ Y\times^R_S Z \ar[r]^{\sim~~~~~~~~~~~~~~} &
      \bE_{Y,Z}[-1]=\bE_T[-1]=\bL_{Y\times^R_S Z} }
  \]
  is commutative, where $E_{X,Z}=E_W=\frac{T_S}{T_X+T_Z}$ and
  $E_{Y,Z}=E_T=\frac{T_S}{T_Y+T_Z}$ are the excess bundles supported on
  \[ W = X\times_S Z \mbox{ and } T = Y \times_S Z.\]
  The horizontal isomorphisms are defined in~\cite{ACH} and will be
  explained in the proof below.  The map
  $\bE_{X,Z}[-1]\rightarrow \bE_{Y,Z}[-1]$ is induced by
  the obvious map of vector bundles.
\end{Proposition}

\paragraph
Before we begin the proof we note that the setup of the above
proposition gives rise to the following diagram of spaces, where $W$
and $T$ are the underived fiber products,
\[
  \xymatrix{ X\ar[dd] & & X\times_S^R Z\ar[dd]\ar[ll] & \\ &
    W=X\times_SZ\ar[ul]\ar[dr]\ar[dl]\ar[ur]\ar[dd] & &\\ Y\ar[dd]& &
    Y\times_S^R Z\ar[dd]\ar[ll]& \\ & T=Y\times_SZ
    \ar[ul]\ar[dr]\ar[dl]\ar[ur] & & \\ S & & Z.\ar[ll] & }
\]

\begin{Proof}
The compatibility of first order splittings implies that the following
two short exact sequences and their splittings are compatible
\[
\xymatrix{ 0\ar[r] & N_{W/Z}\ar[r]\ar[d]
&N_{X/S}|_W\ar@{-->}@/_1.pc/[l]\ar[r]\ar[d] &E_{W}\ar[r]\ar[d] & 0 \\
0\ar[r] & N_{T/Z}|_W\ar[r] & N_{Y/S}|_W\ar@{-->}@/_1.0pc/[l]\ar[r]
&E_{T}|_W\ar[r] & 0.  }
\]
The two isomorphisms $\bE_W[-1]\cong X\times^R_SZ$, and
$\bE_T[-1]\cong Y\times^R_SZ$ are defined using the two
splittings of short exact sequences above.  The three horizontal maps
on the left of the diagram below are the splittings of the short exact
sequences.  The composition of horizontal maps below are the desired
isomorphisms $\bE_W[-1]\cong X\times^R_SZ$ and
$\bE_T[-1]\cong Y\times^R_SZ$,
\[
\xymatrix{ \bE_W[-1]\ar@{-->}[r]\ar[d] &
\bN_{X/S}[-1]|_W\ar[r]^{\sim~~~~~}\ar[d] & X\times_S^R
X|_W=X\times_S^RW\ar[r]\ar[d] & X\times_S^RZ\ar[d] \\
\bE_T[-1]|_W\ar@{-->}[r]\ar[d] &
\bN_{Y/S}[-1]|_W\ar[r]^{\sim~~~~~}\ar[d] & Y\times_S^R
Y|_W=Y\times_S^RW\ar[r]\ar[d] & Y\times_S^RZ\ar[d] & \\
\bE_T[-1]\ar@{-->}[r] &
\bN_{Y/S}[-1]|_T\ar[r]^{\sim~~~~~} & Y\times_S^R
Y|_T=Y\times_S^RT\ar[r] & Y\times^R_S Z. }
\]

We only need to prove the commutativity of the isomorphisms in the
middle
\[
\xymatrix{ \bN_{X/S}[-1]|_W\ar[r]^{\sim~~~~~}\ar[d] &
X\times_S^R X|_W=X\times_S^RW\ar[d]\\
\bN_{Y/S}[-1]|_W\ar[r]^{\sim~~~~~} & Y\times_S^R
Y|_W=Y\times_S^RW.  }
\]
because all the others are commutative.  We can restrict everything to
$X$ first, and then restrict to $W$.  Therefore, it suffices to show
the commutativity of
\[
\xymatrix{ \bN_{X/S}[-1]\ar[r]^{\sim}\ar[d] & X\times_S^R
X\ar[d]\\ \bN_{Y/S}[-1]|_X\ar[r]^{\sim~~~} & Y\times_S^R
Y|_X=Y\times_S^RX. & }
\]
This is Theorem A in \cite{H}.
\end{Proof}

\begin{Proof}[Proof of Proposition \ref{prop assoc 2}.]
For simplicity denote the space $X^g$, $X^h$, $X^k$, and $X$ in
Proposition \ref{prop assoc 2} by $X$, $Y$, $Z$, and $S$.

Because of Proposition \ref{prop assoc 1'} we have the commutativity
of
\[\xymatrix{
X\times^R_SY\ar[r]^{\sim~~~~~~} &
\bE_{X,Y}[-1]=\bL_{X\times^R_SY} &
W\times^R_SZ\ar[r]^{\sim~~~~~~}\ar[d] &
\bE_{W,Z}[-1]=\bL_{W\times^R_SZ}\ar[d]\\
X\times^R_ST\ar[r]^{\sim~~~~~~}\ar[u] &
\bE_{X,T}[-1]=\bL_{X\times^R_ST}\ar[u] &
Y\times^R_SZ\ar[r]^{\sim~~~~~~} &
\bE_{Y,Z}[-1]=\bL_{Y\times^R_SZ}, }
\]
where $T=Y\times_S Z$ and $W=X\times_S Y$.  As a consequence we get the
commutative diagram
\[\mbox{\scriptsize
\xymatrixcolsep{0.2in} \xymatrix{ X\times^R_SY\times^R_SZ
\ar[r]\ar[d]_{id} &
\bE_{X,Y}[-1]\times^R_SZ=\bE_{X,Y}[-1]\times^R_W(W\times^R_S
Z)\ar[r]\ar[d] &
\bE_{X,Y}[-1]\times^R_W\bE_{W,Z}[-1]\ar[r]\ar[d] &
\bE_{X,Y,Z}[-1]\ar[d]_{id}\\ X\times^R_SY\times^R_SZ
\ar[r]\ar[d]_{id} &
\bE_{X,Y}[-1]\times^R_SZ=\bE_{X,Y}[-1]\times^R_Y(Y\times^R_SZ)\ar[r]
&\bE_{X,Y}[-1]\times^R_Y\bE_{Y,Z}[-1]\ar[r]\ar[d]_{id} &
\bE_{X,Y,Z}[-1]\ar[d]_{id} \\
X\times^R_SY\times^R_SZ\ar[r]\ar[u]^{id} &
X\times^R_S\bE_{Y,Z}[-1]=(X\times^R_SY)\times^R_Y\bE_{Y,Z}[-1]\ar[r]
& \bE_{X,Y}[-1]\times_Y\bE_{Y,Z}[-1]\ar[r]\ar[u]^{id} &
\bE_{X,Y,Z}[-1]\ar[u]^{id} \\
X\times^R_SY\times^R_SZ\ar[r]\ar[u]^{id} &
X\times_S^R\bE_{Y,Z}[-1]=(X\times^R_S
T)\times^R_T\bE_{Y,Z}[-1]\ar[r]\ar[u]
&\bE_{X,T}[-1]\times^R_T\bE_{Y,Z}[-1]\ar[r]\ar[u] &
\bE_{X,Y,Z}[-1]\ar[u]^{id} ,} }
\]
where all the arrows are isomorphisms, and $E_{X,Y,Z}$ is the excess
bundle of the triple intersection $X\times^R_SY\times^R_SZ$.

The two rightmost squares of the diagram above commute because there
are natural isomorphisms $E_{X,Y,Z}\cong E_{X,Y}|_U\oplus E_{W,Z}\cong
E_{X,T}\oplus E_{Y,Z}|_U$, where $U=X\cap Y\cap
Z=X\times_SY\times_SZ$.  The commutativity of the outer square of the
diagram above is the one that we needed to prove in Proposition
\ref{prop assoc 2}.
\end{Proof}

\section{A possible simplification}

This section is more speculative. After we rewrite our product in
Definition~\ref{def product} in more concrete terms, we propose a way
to simplify the formulas for Calabi-Yau global quotient orbifolds.
The simplification is motivated by the definition of Chen-Ruan
orbifold cohomology, so we need to first review this definition before
proceeding.  We then provide several examples comparing the simplified
product with the Chen-Ruan orbifold cohomology, via homological mirror
symmetry.  We end the paper by stating a number of questions that
remain open for future research.

\paragraph
\label{subsec:prod1}
In Section~4 we computed the structure complexes of
$\bL_{\tXgch}$ and $\bL_{\tXg}\times^R_X\bL_{\tXh}$, so we can write
the two-step map
\[\bD(\bL_{\tXg})\otimes\bD(\bL_{\tXh}) \ra
  \bD(\bL_{\tXg}\times^R_X\bL_{\tXh})\cong\bD(\bL_{\tXgch}) \ra
  \bD(\bL_\tXgh) \]
in Definition~\ref{def product} in more concrete terms.  Explicitly,
for $g\in G$ define
\[ \HT^{(p,q)}(X;g) = H^{p-c_g}(X^g, \wedge^q
  T_{X^g}\otimes\omega_g). \]
Remark that this is {\em not} the same bigrading as the one in the
introduction. 

The two composition of the maps  above can then be written as a direct
sum over $p, p', q, q'$ of maps
\[  \HT^{(p,q)}(X;g) \otimes \HT^{(p',q')} (X;h) \ra
                             \bigoplus_{i=0}^{\rank E}\HT^{(p+p'-i, q+q'+i)} (X; gh)
\]
factoring through the middle term (coming from $\bD(\bL_\tXgch)$)
\[ \bigoplus_{i=0}^{\rank  E} H^{p+p'-c_{g,h}-i}( X^{g,h},
            \wedge^q T_{X^g}|_{X^{g,h}}\otimes\wedge^{q'}T_{X^h}|_{X^{g,h}}
             \otimes\wedge^iE \otimes\omega_{g,h}). \]
Here $E$ is the excess bundle for the intersection of $X^g$ and
$X^h$ in $X$.  Note in particular that if $G$ is trivial our definition
recovers the classical product on polyvector fields.

Observe that the $\HT^{(p,q)}$ notation does {\em not}
give a bigrading -- {\em a priori} all the maps above, for $0\leq
i\leq \rank E$ could be non-zero.  The simplification we propose, for
the Calabi-Yau case, is to leave only one of these maps, for a
specific $i$.  (Conjecturally, all the other maps would be zero anyway.)

\paragraph{\bf The formulas for Chen-Ruan orbifold cohomology.}
We now discuss some preparations for the motivation for the
simplification of the product we defined.  The idea is to draw
inspiration from mirror symmetry, and to regard, in the Calabi-Yau
case, Chen-Ruan orbifold cohomology as the mirror of orbifold
Hochschild cohomology.

Let $X$ be a complex manifold endowed with the action of a finite
group $G$.  Chen and Ruan~\cite{CR} defined a version of the singular
cohomology {\em ring} for the orbifold $[X/G]$.  Fantechi and
G\"{o}ttsche~\cite{FG} wrote down the formula for the product
explicitly as follows.  They first constructed an associative product
on
\[ H^*_\orb(X;G) = \bigoplus_{g\in G}H^{*-2\iota(g)}(X^g,\C) \]
which maps $\alpha_g\in H^{*-2\iota(g)}(X^g,\C)$ and
$\beta_h\in H^{*-2\iota(h)}(X^h,\C)$ to
\[ (\alpha_g, \beta_h) \mapsto
  i^{gh}_{g,h*}(\alpha_g|_{X^{g,h}}\cdot\beta_h|_{X^{g,h}}\cdot\gamma_{g,h}). \]
Here $\gamma_{g,h}$ is the top Chern class of a certain twist bundle
whose rank is $\iota(g)+\iota(h)-\iota(gh)-\codim(X^{g,h},X^{gh})$,
where $\iota(g)$ is the so-called {\em age} of $g$, see~\cite{FG}.

The Chen-Ruan orbifold singular cohomology ring is obtained by taking
$G$-invariants: 
\[ H^*_\orb([X/G]) = H^*_\orb(X;G)^G.  \]

Note that the above ring is {\em bigraded} with respect to the orbifold
Hodge decomposition~\cite{ALR}
\[H^{n-2\iota(g)}(X^g,\C) =
  \bigoplus_{p+q=n}H^{p-\iota(g)}(X,\wedge^{q-\iota(g)}\Omega_{X^g}).\] 

\paragraph
Mirror symmetry associates two graded commutative rings to a
Calabi-Yau space: the A- and the B-model state spaces, which are
interchanged by the mirror operation.  When the target space is a
compact Calabi-Yau manifold $X$, the A-space is $H^*(X, \C)$, while
the B-space is $\HH^*(X)$.  When it is an orbifold $[X/G]$, these
spaces are naturally the Chen-Ruan orbifold cohomology and the
orbifold Hochschild cohomology rings of $[X/G]$.  

Since the product in the A-model preserves the $p,q$ bidegree, the
yoga of mirror symmetry suggests that the product in the B-model
should also preserve {\em some} bidegree for Calabi-Yau orbifolds.

The proofs of the following two lemmas are left as exercises to the reader.

\begin{Lemma}
  \label{lem iso}
  There is a natural isomorphism
  \[\omega_{g}|_{X^{g,h}}[-c_g]\otimes\omega_{h}|_{X^{g,h}}[-c_h] \iso
    \wedge^r  E[r]\otimes\omega_{g,h}[-c_{g,h}],\]
  where $r$ is the rank of the excess bundle $E$.
\end{Lemma}
\medskip

\begin{Lemma}
  \label{lem contract}
  The bundle $T_{X^g}|_{X^{g,h}}$ decomposes naturally into a direct
  sum as $T_{X^{g,h}}\oplus N_{X^{g,h}/X^g}$, and similarly for
  $T_{X^h}|_{X^{g,h}}$. 

  The class $\gamma_{g,h}\in H^k(X^{g,h},\wedge^k\Omega_{X^{g,h}})$ in
  Fantechi and G\"{o}ttsche's paper~\cite{FG} acts naturally
  on
  \[ \bigoplus_{p,q} H^{p}(X^{g,h},\wedge^{q}T_{X^g}|_{X^{g,h}}\otimes\wedge^{q'}T_{X^h}|_{X^{g,h}}
    \otimes\wedge^r E^{\chk}\otimes\omega_{g,h}), \]
  where
  \[ k=\iota(g)+\iota(h)-\iota(gh)-\codim(X^{g,h},X^{gh}) \]
  and the action is given by the contraction of $\Omega_{X^{g,h}}$
  with $T_{X^{g,h}}$.
\end{Lemma}
\medskip

\paragraph
\label{def simplify} 
We are now ready to give a new construction for an operation on
$\HT^*(X; G)$ which mimics more closely the Fantechi-G\"ottsche
product~\cite{FG}.  Define the bigraded piece $\HT^{p,q}(X;G)$ of
bidegree $p,q$ of $\HT(X;G)$ by
\[ \HT^{p,q}(X;G) =
  H^{p-\iota(g)}(X,\wedge^{q+\iota(g)-c_g}T_{X^g}\otimes\omega_g). \]
The product will be bigraded, being given by maps
\[ \HT^{p,q}(X;G) \otimes \HT^{p',q'}(X;G) \ra \HT^{p+p',
    q+q'}(X;G). \]
Note that unlike the product in~(\ref{subsec:prod1}), only one of the
maps there is non-zero.  We conjecture that in Calabi-Yau situations,
the two products agree -- in other words, all the maps
in~(\ref{subsec:prod1}) which do not preserve the bigrading are zero.
This is the case in all the examples we study below.

The new product is defined as the following composition:
 \begin{align*}
H^{p}  (X^g, & \wedge^q T_{X^g}\otimes\omega_g[-c_g])\otimes
             H^{p'}(X^h,\wedge^{q'}T_{X^h}\otimes\omega_h[-c_h]) \\
&\ra  H^{p+p'}(X^{g,h},\wedge^q
T_{X^g}|_{X^{g,h}}\otimes\omega_{g}|_{X^{g,h}}[-c_g] \otimes
   \wedge^{q'}T_{X^h}|_{X^{g,h}}\otimes\omega_h|_{X^{g,h}}[-c_h])
   \\ 
& \cong H^{p+p'-r}(X^{g,h},\wedge^q T_{X^g}|_{X^{g,h}}\otimes
\wedge^{q'}T_{X^h}|_{X^{g,h}}\otimes\omega_{g,h}[-c_{g,h}]\otimes\wedge^rE)
   \\
&  \ra\bigoplus_{i+j=k}H^{p+p'-r+k}(X^{g,h},\wedge^{q-i}
T_{X^g}|_{X^{g,h}}\otimes
\wedge^{q'-j}T_{X^h}|_{X^{g,h}}\otimes\omega_{g,h}[-c_{g,h}]\otimes\wedge^rE)
   \\
& \ra H^{p+p'-r+k}(X^{gh},\wedge^{q+q'+r-k}T_{X^{gh}}\otimes\omega_{gh}[-c_{gh}]).
\end{align*}
The first arrow is the naive restriction from $X^g$ and $X^h$ to
$X^{g,h}$.  The isomorphisms in the middle are due to Lemma~\ref{lem
iso}.  The last arrow is the map ${\bL_{m}}_*$ in Definition~\ref{def
product}.  The second arrow in the middle involving $k$ is the action
of $\gamma_{g,h}$ in Lemma~\ref{lem contract}.  One does indeed verify
that this map respects the bigrading defined above.

\paragraph{\bf Examples.}
For a first example consider an abelian surface $A$ endowed with the
action of $\Z/2\Z$, acting by negation in the group law of $A$.  The
tangent bundle of $A$ is trivial, so there are no Duflo correction
terms.  The mirror of the orbifold $[A/G]$ is expected to be $[A/G]$
itself in this case.   This suggests that the
product we defined should match with the one on the orbifold
cohomology of $[A/G]$, so we expect to find isomorphisms
\[\HT^{*}([A/G]) \iso \HH^{*}([A/G])\cong H^*_{\orb}([A/G],\C),\]

It is known~\cite{FG} that in this case the classes $\gamma_{g,h}$ are
trivial. Write $G = \{e,\tau\}$ where $e$ is the identity element.
Then we have 
\[\HH^*([A/G])=\left (\HH^*(A,e)\oplus\HH^*(A,\tau)\right)^G=
  \HH^*(A,e)\oplus\HH^*(A,\tau)^\tau,\]
where for $g\in G$ the notation $\HH^*(A,g)$ was explained in
Section~2.  The space $\HH^*(A,e)$ is the Hochschild cohomology of
$A$, and its product is well-understood from the Kontsevich and
Calaque-Van den Bergh theorem.  The only non-trivial product we need
to understand is
\[ \HH^*(A,\tau)\otimes\HH^*(A,\tau)\rightarrow\HH^*(A,e). \]

Note that the space
\[ \HH^*(A,\tau)=H^{0}(A^{\tau},\wedge^0T_{A^\tau}\otimes\omega_\tau)=
  H^0(A^\tau,\C), \]
is a $16$-dimensional vector space in cohomological degree $2$.
It is of bidegree $(1,1)$ under the new bigrading we defined
in~(\ref{def simplify}).  By the definition of our product, it 
is also clear that the product of two $(1,1)$-form gives a
$(2,2)$-form which lands in $H^2(A,\wedge^2T_{A})$.  This matches
perfectly with the product on orbifold cohomology~\cite{FG}.

\paragraph
For another example, consider a holomorphic symplectic orbifold
$[X/G]$.  Again, the mirror of $[X/G]$ is expected to be $[X/G]$, so
we expect to get
\[ \HT^*([X/G])\cong\HH^*(X/G)\cong H^*_{\orb}([X/G],\C).\]

The right hand side decomposes into
\[ H^{*-2\iota(g)}(X^g,\C)=\bigoplus_{p+q=*}
  H^{p-\iota(g)}(X^g,\wedge^{q-\iota(g)}\Omega_{X^g})\]
by the Hodge decomposition.  The left hand side is \[\bigoplus_{g\in
G}H^{p-\iota(g)}(X^g,\wedge^{q+\iota(g)-c_g}T_{X^g}\otimes\omega_g).\]

Moreover, $\omega_g$ is trivial and $\Omega_{X^g}\cong T_{X^g}$
because of the holomorphic symplectic condition.  There is a canonical
identification between the two sides as vector spaces, and we believe
the two products should agree.  The bigradings of the two sides match
completely because $2\iota(g)=c_g$.

One important example we have in mind is when $X=K^n$ consists of $n$
copies of a K3 surface $K$ and $G=\Sigma_n$ is the symmetric group
acting on $K^n$ by permutation.  The group is not abelian in this
case, but the constructions and results in Sections~3 and~4 still
work because one can check directly that Proposition \ref{Prop
  formality'} holds in this situation.  The key point is that in this
situation all the tangent bundles and normal bundles involved are
copies of direct sums of $T_K$, so the short exact sequence in the
proof of Proposition~\ref{Prop formality'} splits naturally.

\paragraph{\bf Open questions.}

$(1)$ We can not prove that the simplified product in Definition
\ref{def simplify} agrees with the one in Definition~\ref{def
product}.  We conjecture that they agree under the Calabi-Yau assumption.

$(2)$ For any Calabi-Yau orbifold, we believe that our product on
orbifold polyvector fields should match with the Chen-Ruan orbifold 
cohomology of the mirror.

$(3)$ This is our main Conjecture A.  For any orbifold $[X/G]$ with an
abelian group action, we believe that Kontsevich's Theorem holds,
i.e., the orbifold Hochschild cohomology should be isomorphic to the orbifold
polyvector fields.   More precisely, we conjecture that the diagram
\[
\xymatrix{ \bD(\bL_{\tXg})\otimes\bD(\bL_{\tXh})\ar[d]
\ar[rrrr]^{HKR\circ(\sqrt{\Td(T_{X^g})}\lrcorner-)\otimes
HKR\circ(\sqrt{\Td(T_{X^h})}\lrcorner-)} & & & &
\bD(\tXg)\otimes\bD(\tXh)\ar[d]\\
\bD(\bL_{\tXg\times^R_X\tXh})\cong\bD(\bL_{\tXg}\times^R_X\bL_{\tXh})
\ar[d]^{{\bL_m}_*}\ar[rrrr]^{HKR\circ(\sqrt{\Td(T_{X^{g,h}})}\lrcorner-)}
& & & & \bD(\tXg\times^R_X\tXh)=\bD(\tXgch)\ar[d]^{m_*}\\
\bD(\bL_{\tXgh})
\ar[rrrr]^{HKR\circ(\sqrt{\Td(T_{X^{gh}})}\lrcorner-)} & & & &
\bD(\tXgh) } \]
is commutative, where the horizontal maps are isomorphisms.  All the
HKR maps that appear in the horizontal isomorphisms are the formality
isomorphisms in Sections~2--4.  They generalize the classical HKR
isomorphism as explained in $\cite{ACH}$.  As mentioned at the very
beginning of this paper, HKR can not be an isomorphism of rings, so we
need to add the Duflo correction term in the horizontal isomorphisms.

\end{document}